\documentclass{scrartcl}
\usepackage[ngerman, english]{babel}
\usepackage{amsmath}
\usepackage{yhmath}
\usepackage{amssymb}
\usepackage{amsfonts}
\usepackage[all]{xy}
\usepackage[amsmath,thmmarks,thref]{ntheorem}              
\usepackage{mathtools}    
\usepackage{mathrsfs}     
\usepackage{bbding}
\usepackage{cancel}       
\usepackage{relsize}
\usepackage{paralist}
\usepackage{nicefrac}        
\usepackage{graphicx}
\usepackage{tikz}
\usepackage{subfigure}

\font\matha = matha10

\DeclarePairedDelimiter\abs{\lvert}{\rvert}
\DeclarePairedDelimiter\norm{\lVert}{\rVert}

\providecommand{\N}{\mathbb{N}}
\providecommand{\Z}{\mathbb{Z}}
\providecommand{\R}{\mathbb{R}}

\providecommand{\eins}{\mathbf 1}

\providecommand{\inv}{\bar{~}}                         
\providecommand{\dual}{^{\mathsmaller{\mathsmaller{\wedge}}}}    

\providecommand{\gdw}{\Leftrightarrow}

\providecommand{\eps}{\varepsilon}

\DeclareMathOperator{\supp}{supp}   

\newcommand{\conv}{\mathbin{\mbox{\matha\char006}}}

\theoremstyle{break}                
\theoremheaderfont{\normalfont\bfseries}
\theorembodyfont{\itshape}
\theoremseparator{:}
\theoremnumbering{arabic}
\theoremsymbol{}

\newtheorem{pro}{Proposition}
\newtheorem{satz}[pro]{Theorem}

\theoremstyle{break}           
\theoremheaderfont{\normalfont\bfseries} 
\theorembodyfont{\upshape}
\theoremseparator{:}
\theoremnumbering{arabic}
\theoremsymbol{}

\newtheorem{dfn}[pro]{Definition}

\theoremstyle{nonumberplain}
\theoremheaderfont{\scshape}
\theorembodyfont{\upshape}
\theoremseparator{:}
\theoremnumbering{arabic}
\theoremsymbol{\ensuremath{\square}}
\newtheorem{bew}{Proof}

\theoremstyle{nonumberplain}

\numberwithin{equation}{section}

\setdefaultleftmargin{2.5em}{2.2em}{1.87em}{1.7em}{1em}{1em}
 
\usepackage[utf8]{inputenc}
\usepackage{makeidx}
\usepackage{graphicx}
\usepackage{fancyhdr}
\usepackage{lmodern}
\addtokomafont{disposition}{\rmfamily}
\addtokomafont{section}{\rmfamily}
\addtokomafont{title}{\rmfamily}
\addtokomafont{subsection}{\rmfamily}
\addtokomafont{paragraph}{\rmfamily}

\mathtoolsset{
showonlyrefs,
showmanualtags}
\makeindex
\title{About a theorem of Wiener on the Bessel-Kingman Hypergroup}

\author{Lukas Innig\\\ttfamily lukas.innig@googlemail.com}

\begin{document}


\maketitle

\begin{abstract}
  A theorem of Wiener on the circle group was strengthened and extended by Fournier in \cite{fournier1} to locally compact abelian groups and extended further to the Bessel-Kingman hypergroup with parameter $\alpha=\nicefrac12$ by Bloom/Fournier/Leinert in \cite{leinpaper}. We further extend this theorem to Bessel-Kingman hypergroups with parameter $\alpha > \nicefrac12$.
\end{abstract}

\section{The Bessel-Kingman Hypergroup}

In this paper we will prove a theorem of \emph{Fournier} (\thref{satzfournier}) on the Bessel-Kingman Hypergroup with parameter $\alpha\geq\frac12$. We will use the proof in \cite{leinpaper} for the case $\alpha=\frac12$ as a guideline, which will be altered where necessary. \thref{satzfournier} is based upon a theorem of \emph{Wiener} but treats a more general case. Following \cite{bloomheyer} and \cite{leinpaper} we define:

\begin{dfn}
The Bessel-Kingman Hypergroup with parameter $\alpha$ is defined as $K=(\R_+,\conv_\alpha)$, where 
  \begin{equation}
    \eps_x\conv_\alpha\eps_y(f)=\int_{\abs{x-y}}^{x+y}K_\alpha(x,y,z) f(z)z^{2\alpha+1}\,dz
  \end{equation}
with
\begin{equation}
  K_\alpha(x,y,z)=\frac{\Gamma(\alpha+1)}{\Gamma(\nicefrac12)\Gamma(\alpha+\nicefrac12)2^{2\alpha-1}}\frac{\left[\left(z^2-(x-y)^2\right)\left((x+y)^2-z^2\right)\right]^{\alpha-\nicefrac12}}{(xyz)^{2\alpha}}
\end{equation}
and identity involution ($x^-=x$). The Haar measure $\omega_\alpha(dz)$ is of the form $\omega_\alpha(dz)=z^{2\alpha+1}dz$. The characters are given by $\chi_\lambda(x)\coloneqq j_\alpha(\lambda x),\;x\in\R_+$ where $j_\alpha$ denotes the modified Bessel function of order $\alpha$:
\begin{equation}
  j_\alpha(x)\coloneqq\sum_{k=0}^\infty\frac{(-1)^k\Gamma(\alpha+1)}{2^{2k}k!\Gamma(\alpha+k+1)}x^{2k},\quad \forall x\in \R.
\end{equation}
One has $\chi_0\equiv 1$. Furthermore $K$ is a Pontryagin hypergroup. In fact $K\cong K\dual$, where the isomorphism is given by $\lambda\mapsto\chi_\lambda$. We note that $(\R_+,\conv_\alpha)$ is commutative because $K(x,y,z)=K(y,x,z)$ and $\int_{\abs{x-y}}^{x+y}\cdots=\int_{\abs{y-x}}^{y+x}\cdots$.
\end{dfn}

As a convention, we denote
\begin{itemize}
\item $I_n\coloneqq [n-1,n)$,
\item $C_\Gamma\coloneqq\frac{\Gamma(\alpha+1)}{\Gamma(\nicefrac12)\Gamma(\alpha+\nicefrac12)2^{2\alpha-1}}$,
\item $\omega_n\coloneqq\omega_\alpha(I_n)$.
\end{itemize}

Furthermore let $\alpha\geq\nicefrac12$. The case $-\nicefrac12<\alpha<\nicefrac12$ will not be treated.

\begin{dfn}\label{gewoehnlichenorm}
For a hypergroup $(\R_+,\conv)$ with Haar measure $\omega$, the \emph{discrete amalgam norm} is given by
  \begin{equation}
    \norm f_{p,q}\coloneqq \left(\sum_{n=1}^\infty\omega_n\left(\frac1{\omega_n}\int_{I_n}\abs f^p\;d\omega\right)^{\nicefrac q p}\right)^{\nicefrac 1 q}.
  \end{equation}
In the case $p$ or $q$ equal to $\infty$ we set by convention
\begin{align}
  \norm f_{\infty,q}&\coloneqq \left(\sum_{n=1}^\infty\omega_n\sup_{x\in I_n}\abs{f(x)}^q\right)^{\nicefrac1q},\\
  \norm f_{p,\infty}&\coloneqq \sup_{n\in\N}\left(\frac1{\omega_n}\int_{I_n}\abs f^p\;d\omega\right)^{\nicefrac1p}\quad\text{and }\\
  \norm f_{\infty, \infty}&\coloneqq \sup_{n \in \N}\left(\sup_{x \in I_n}\abs{f(x)}\right) = \norm f_\infty.
\end{align}
The function spaces $\{f\text{ measurable}\mid\norm f_{p,q}<\infty\}$ will be denoted as $(L^p,\ell^q)(\R_+,\conv)$. For these spaces the following properties hold:

\begin{align}
  \norm f_{p_1,q}&\leq\norm f_{p_2,q},&&\quad\text{if } p_1\leq p_2.\\
  \norm f_{p,q_1}&\leq C\norm f_{p,q_2},&&\quad\text{if } q_1\geq q_2.\\
\intertext{particularly, it holds for  $p_1\leq p_2$ and $q_1\geq q_2$}\\
(L^{p_2},\ell^{q_2})(\R_+,\conv)&\subset(L^{p_1},\ell^{q_1})(\R_+,\conv)\\
\shortintertext{and}
(L^p,\ell^q)(\R_+,\conv)&\subset L^p(\R_+,\conv)\cap L^q(\R_+,\conv)&&\quad\text{for } p\geq q,\\
L^p(\R_+,\conv)\cup L^q(\R_+,\conv)&\subset (L^p,\ell^q)(\R_+,\conv)&&\quad\text{for } p\leq q.
\end{align}
\end{dfn}

\begin{dfn}\label{intrinsischenorm}
Because we are situated on a hypergroup, we can also form amalgam spaces by shifting the unit interval $I_1$ using the left-translation $\tau_y$ defined as
\begin{equation}
\tau_yf(x) = f(y\conv_\alpha x) = \int_{\abs{x-y}}^{x+y}K_\alpha(x,y,z) f(z)z^{2\alpha+1}\,dz\;.
\end{equation}
For the Bessel-Kingman hypergroup $(\R_+,\conv_\alpha)$ the \emph{continuous $(p,\infty)$-amalgam norm} is given by
  \begin{equation}
    \sup_{y\in\R_+}\left(\int\abs f^p\tau_y\eins_{[0,1)}\;d\omega\right)^{\nicefrac1p}.
  \end{equation}
\end{dfn}

Now we are able to state our more general version of the theorem of Fournier \cite[Theorem 3.1]{fournier1}. The proof of this theorem is the main part of this paper. As said before it is closely related to a theorem of Wiener which implies that an integrable function with non-negative Fourier transform on the unit circle which is square integrable on a unit neighborhood is already square integrable on the whole circle. In \cite{fournier1} this theorem was extended to $\R$ and to LCA-Groups.

\begin{satz}[Theorem of Fournier]\label{satzfournier}
  For $f\in L^1(\R_+,\conv_\alpha)$ with $\hat f\geq0$ the following statements are equivalent:
  \begin{enumerate}
  \item $f$ is square integrable on a neighborhood of $0$.
  \item $\hat f\in(L^1,\ell^2)(\R_+,\conv_\alpha)$.
  \item $f\in(L^2,\ell^\infty)(\R_+,\conv_\alpha)$.
  \end{enumerate}
\end{satz}
To prove this we follow \cite{leinpaper}. So we have to check the following properties of $(\R_+,\conv_\alpha)$ with $\alpha\geq\nicefrac12$.

\paragraph{Equivalence of the discrete and the continuous amalgam norms.} 
We show that the norms defined in \thref{gewoehnlichenorm} and \thref{intrinsischenorm} are equivalent.

\paragraph{Uniform boundedness of the translation operator on $(L^p,\ell^q)(\R_+,\conv_\alpha)$.}

We show that
\begin{equation}
  \norm{\tau_yf}_{p,q}\leq C\cdot\norm f_{p,q},\quad \forall y\in\R_+
\end{equation}
with a constant $C$ independent of $y$. For this we will first prove the uniform boundedness on $(L^\infty,\ell^1)$. Then by invoking duality and interpolation arguments we get the general case.

\paragraph{Properties of the convolution.}

For the convolution a generalized version of the Young inequality on amalgams must hold:
\begin{satz}[Young inequality on amalgams]
  For $f_1\in(L^{p_1},\ell^{q_1})$ and $f_2\in(L^{p_2},\ell^{q_2})$, where
  \begin{equation}
    \left(\frac1p,\frac1q\right)=\left(\frac1{p_1},\frac1{q_1}\right)+\left(\frac1{p_2},\frac1{q_2}\right)-(1,1)
  \end{equation}
we have $f_1\conv f_2\in(L^p,\ell^q)$ and
\begin{equation}
  \norm{f_1\conv f_2}_{p,q}\leq C\cdot\norm{f_1}_{p_1,q_1}\cdot\norm{f_2}_{p_2,q_2}.
\end{equation}
\end{satz}

\paragraph{Properties of the Fourier transformation.}

For the Fourier transform of a function $f$ we need a more generalized form of the Hausdorff-Young theorem: 
\begin{satz}[Hausdorff-Young theorem on amalgams]\label{satzhy}
  For $f\in(L^p,\ell^q)(\R_+,\conv_\alpha)$ with $1\leq p,q\leq 2$ it holds, that $\hat f\in(L^{q'},\ell^{p'})(\R_+,\conv_\alpha)$.
\end{satz}
Here too, we will consider the special cases $(p,q)\in\{(1,1),\,(1,2),\,(2,1),\,(2,2)\}$. The result then follows with interpolation arguments.

\subsection{Equivalence of the discrete and the continuous amalgam norms}

\begin{pro}
   We have
   \begin{equation}
     \norm f_{p,\infty}\leq C\sup_{y\in\R_+}\left(\int\abs f^p\tau_y\eins_{[0,1]}\;d\omega_\alpha\right)^{\nicefrac1p}.\label{contnormdef}
   \end{equation}
\end{pro}
\begin{bew}
Let us concentrate on the term $\tau_y\eins_{[0,1]}$:
  \begin{equation}
    \tau_y\eins_{[0,1]}(x)=\frac{C_\Gamma}{(xy)^{2\alpha}}\int_{\abs{x-y}}^{x+y}\eins_{[0,1]}(z)z\left[\left(z^2-(x-y)^2\right)\left((x+y)^2-z^2\right)\right]^{\alpha-\nicefrac12}\;dz.
  \end{equation}
Like in \cite{leinpaper}, it will be sufficient to look at the supremum in \eqref{contnormdef} taken only over all $y$ of the form $y=n+\nicefrac12,\;1\leq n\in\Z$. We therefore consider $\tau_{n+\nicefrac12}\eins_{[0,1]}(x)$ for $x\in I_{n+1}$. The domain of the integration is
 \begin{align}
    [\abs{x-y},x+y]\cap[0,1]&=[\abs{x-n-\nicefrac12},x+n+\nicefrac12]\cap[0,1]\\
    &=[\abs*{x-n-\nicefrac12},1]\supset[\nicefrac12,1] \qquad\text{for
    }x\in I_{n+1}.
  \end{align}
By using  $x^{2\alpha}=x\cdot (x^2)^{\alpha-\nicefrac12}$ together with  $\abs{x-n-\nicefrac12}\leq \nicefrac12$ we get
  \begin{align}
    \MoveEqLeft\frac{C_\Gamma}{(n+\nicefrac12)^{2\alpha}x^{2\alpha}}\int_{\abs{x-n-\nicefrac12}}^1z\left[\left(z^2-(x-n-\nicefrac12)^2\right)\left((x+n+\nicefrac12)^2-z^2\right)\right]^{\alpha-\nicefrac12}\;dz\\
    & \geq\frac{C_\Gamma}{(n+\nicefrac12)^{2\alpha}}\frac 1 x\int_{\nicefrac12}^1z\Bigg[\left(z^2-\frac14\right)\underbrace{\left(\frac{(x^2+2x(n+\nicefrac12)+(n+\nicefrac12)^2-z^2}{x^2}\right)}_{=(\ast)}\Bigg]^{\alpha-\nicefrac12}\;dz.\\
    (\ast)& =1+\frac{2(n+\nicefrac12)}x+\frac{(n+\nicefrac12)^2-z^2}{x^2}\text{
     and $\frac1x$ are decreasing in $x$.}
  \end{align}

  Therefore it follows that $\tau_{n+\nicefrac12}\eins_{[0,1]}(x)\geq\tau_{n+\nicefrac12}\eins_{[0,1]}(n+1)\;\forall x\in I_{n+1}$. Further we will show that
  \begin{equation}
\tau_{n+\nicefrac12}\eins_{[0,1]}(n+1) \geq C\cdot \frac{1}{n^{2\alpha+1}}.
\end{equation}

  \begin{align}
    \MoveEqLeft\tau_{n+\nicefrac12}\eins_{[0,1]}(n+1)\\
    & =\frac{C_\Gamma}{((n+1)(n+\nicefrac12))^{2\alpha}}\\
    & \qquad\cdot \int_{\nicefrac12}^1z\left[\left(z^2-\frac14\right)\left(\left(2n+\frac32\right)^2-z^2\right)\right]^{\alpha-\nicefrac12}\;dz\\
\shortintertext{(sort by powers of $n$)}
    & =\frac{C_\Gamma}{((n+1)(n+\nicefrac12))^{2\alpha}}\\
    & \qquad \cdot\int_{\nicefrac12}^1z\left[n^2(4z^2-1)+n\left(6z^2-\frac32\right)+\left(\frac{10}4z^2-z^4-\frac9{16}\right)\right]^{\alpha-\nicefrac12}\;dz\\
\shortintertext{(place $n^2$ outside the brackets)}
    & =\frac{C_\Gamma\cdot n^{2\alpha-1}}{((n+1)(n+\nicefrac12))^{2\alpha}}\\
    & \qquad\cdot\int_{\nicefrac12}^1z\left[(4z^2-1)+n^{-1}\left(6z^2-\frac32\right)+n^{-2}\left(\frac{10}4z^2-z^4-\frac9{16}\right)\right]^{\alpha-\nicefrac12}\;dz\\
    & \geq \frac{C_\Gamma\cdot n^{2\alpha-1}}{((n+1)(n+\nicefrac12))^{2\alpha}}\underbrace{\int_{\nicefrac12}^1z(4z^2-1)^{\alpha-\nicefrac12}\;dz}_{\text{constant}}=C'\cdot \frac{n^{2\alpha-1}}{((n+1)(n+\nicefrac12))^{2\alpha}}\\
    & \geq C'\cdot
    \frac{n^{2\alpha-1}}{((n+n)(n+n))^{2\alpha}}=C'\cdot
    \frac{n^{2\alpha-1}}{(4n^2)^{2\alpha}} = C\cdot
    \frac{1}{n^{2\alpha+1}}.
  \end{align}
Above we used that both $f_1(z)\coloneqq\left(6z^2-\frac32\right)$ and $f_2(z)\coloneqq\left(\frac{10}4z^2-z^4-\frac9{16}\right)$ are greater than $0$ on the interval $[\nicefrac12,1]$.\\
\begin{align}
    f_1(z)\geq 0 & \gdw 6z^2-\frac32\geq0\gdw z^2\geq\frac14\gdw\abs z\geq \frac12\\
    f_2(z)\geq 0 & \gdw \frac{10}4z^2-z^4-\frac9{16}\geq 0\\
    & \gdw -\left(z^2-\frac 54\right)^2+1\geq0\gdw1\geq \abs{z^2-\frac 54}\\
    & \gdw\frac32\geq\abs z\geq \frac12.
  \end{align}
   \begin{figure}
     \centering
     \includegraphics{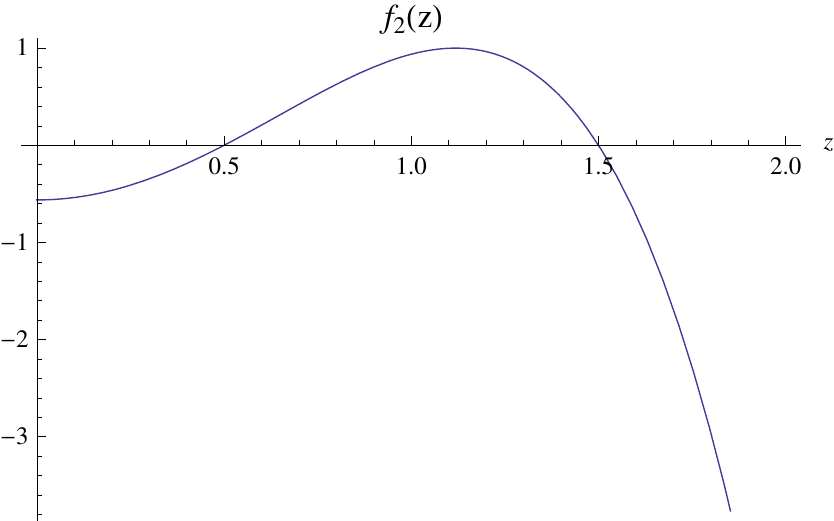}
     \caption{$f_2(z)$}
   \end{figure}
Further we have
  \begin{equation}
    \omega_{n+1}=\int_n^{n+1}z^{2\alpha+1}\;dz\geq n^{2\alpha+1}.
  \end{equation}
  So
  \begin{align}
    \tau_{n+\nicefrac12}\eins_{[0,1]}(n+1) & \geq C\cdot \frac{1}{n^{2\alpha+1}}\\
    & \geq C\cdot\frac1{\omega_{n+1}}\quad\text{for $n\geq1$.}
  \end{align}
Now we can summarize.
  \begin{align*}
    \sup_{y\in[1,\infty)}\left(\int\abs f^p\tau_y\eins_{[0,1]}\;d\omega_\alpha\right)^{\nicefrac1p} & \geq \sup_{n\geq1}\left(\int_{n-\nicefrac12}^{n+\nicefrac32}\abs f^p\tau_{n+\nicefrac12}\eins_{[0,1]}\;d\omega_\alpha\right)^{\nicefrac1p}\\
    & \geq \sup_{n\geq 1}\left(\int_{I_{n+1}}\abs f^p\tau_{n+\nicefrac12}\eins_{[0,1]}\;d\omega_\alpha\right)^{\nicefrac1p}\\
    & \geq\sup_{n\geq 1}\left(\tau_{n+\nicefrac12}\eins_{[0,1]}(n+1)\int_{I_{n+1}}\abs f^p\;d\omega_\alpha\right)^{\nicefrac1p}\\
    & \geq\sup_{n\geq 1}\left(C\frac1{\omega_{n+1}}\int_{I_{n+1}}\abs f^p\;d\omega_\alpha\right)^{\nicefrac1p},\\
    \shortintertext{moreover for $x,y\in I_1$ it obviously holds that $\sup_{y\in
    I_1}\tau_y\eins_{[0,1]}(x)\geq \tau_0\eins_{[0,1]}(x)=1$. So}
    \sup_{y\in\R_+}\left(\int\abs f^p\tau_y\eins_{[0,1]}\;d\omega_\alpha\right)^{\nicefrac1p} & \geq C\cdot\norm{f}_{p,\infty}.
  \end{align*}
\end{bew}

\begin{pro}
  It holds that
  \begin{equation}
    \norm{f}_{p,\infty}\geq C\sup_{y\in\R_+}\left(\int\abs f\tau_y\eins_{[0,1]}\;d\omega_\alpha\right)^{\nicefrac1p}.
  \end{equation}
\end{pro}
\begin{bew}
  \begin{enumerate}
\item Let $y\in [0,1)$. Then one has
    \begin{equation}
      \tau_y\eins_{[0,1]}(x)
      \begin{cases}
        =1, & \text{if } x+y\leq 1,\\
        \leq 1, & \text{if } \abs{x-y}\leq 1\text{ and } x+y\geq 1,\\
        =0, & \text{if } \abs{x-y}\geq 1.
      \end{cases}
    \end{equation}
Note that  $\tau_y\eins_{[0,1]}(x)=\eins_{[0,1]}(x\conv y)=\eps_x\conv\eps_y(\eins_{[0,1]})$. Therefore if $\supp(\eps_x\conv\eps_y)=[\abs{x-y},x+y]\subset [0,1]$, then $\eps_x\conv\eps_y(\eins_{[0,1]})=1$ because $\eps_x\conv\eps_y\in M^1(K)$.\\
It follows now that $\tau_y\eins_{[0,1]}\leq\eins_{[0,2)}$ and therefore
    \begin{align}
      \int\abs f^p\tau_y\eins_{[0,1]}\;d\omega_\alpha & \leq \int\abs f^p\eins_{[0,1)}\;d\omega_\alpha+\int\abs f^p\eins_{[1,2)}\;d\omega_\alpha\\
      & \leq \int\frac1{\omega_1}\eins_{[0,1)}\abs f^p\;d\omega_\alpha+2^{2\alpha+1}\int\frac1{\omega_2}\eins_{[1,2)}\abs f^p\;d\omega_\alpha\\
\intertext{because $\omega_n=\int_{n-1}^{n}z^{2\alpha+1}\;dz\leq n^{2\alpha+1}$, hence $1\leq\frac{n^{2\alpha+1}}{\omega_n}$.}
      \text{So } \left(\int\abs f^p\tau_y\eins_{[0,1]}\;d\omega_\alpha\right)^{\nicefrac1p}&\leq\left(\int_{I_1}\frac1{\omega_1}\abs f^p\;d\omega_\alpha\right)^{\nicefrac1p}+2^{\frac{2\alpha+1}p}\left(\int_{I_2}\frac1{\omega_2}\abs f^p\;d\omega_\alpha\right)^{\nicefrac1p}\\
      &\leq (1+2^{\frac{2\alpha+1}p})\norm f_{p,\infty}.
    \end{align}
.
  \item If $y\in [1,2)$ then it follows, analogously to the calculation above, that $\tau_y\eins_{[0,1]}\leq\eins_{[0,3)}$.
    \begin{align}
      \int\abs f^p\tau_y\eins_{[0,1]}\;d\omega_\alpha & \leq \int\abs f^p\eins_{[0,1)}\;d\omega_\alpha+\int\abs f^p\eins_{[1,2)}\;d\omega_\alpha+\int\abs f^p\eins_{[2,3)}\;d\omega_\alpha\\
      \begin{split}
       & \leq\int\frac1{\omega_1}\eins_{[0,1)}\abs f^p\;d\omega_\alpha+2^{2\alpha+1}\int\frac1{\omega_2}\eins_{[1,2)}\abs f^p\;d\omega_\alpha\\
       &  \qquad +3^{2\alpha+1}\int\frac1{\omega_3}\eins_{[2,3)}\abs f^p\;d\omega_\alpha,
      \end{split}\\
\intertext{therefore we have}
      \begin{split}
       \left(\int\abs f^p\tau_y\eins_{[0,1]}\;d\omega_\alpha\right)^{\nicefrac1p} & \leq\left(\int_{I_1}\frac1{\omega_1}\abs f^p\;d\omega_\alpha\right)^{\nicefrac1p}\\
        & \qquad + 2^{\frac{2\alpha+1}p}\left(\int_{I_2}\frac1{\omega_2}\abs f^p\;d\omega_\alpha\right)^{\nicefrac1p}\\
        & \qquad +3^{\frac{2\alpha+1}p}\left(\int_{I_3}\frac1{\omega_3}\abs f^p\;d\omega_\alpha\right)^{\nicefrac1p}
      \end{split}\\
      &\leq (1+2^{\frac{2\alpha+1}p}+3^{\frac{2\alpha+1}p})\norm f_{p,\infty}.
    \end{align}

  \item Let $k\geq 2$. We consider  $y\in I_{k+1}=[k,k+1]$\\
    $\Rightarrow \supp \tau_y\eins_{[0,1]}=[y-1,y+1]\subset[k-1,k+2]$.\\
   Now it follows 
    \begin{align}
      \tau_y\eins_{[0,1]}(x) & =\frac{C_\Gamma}{x^{2\alpha}y^{2\alpha}}\int_{\abs{x-y}}^1z\left[(z^2-(x-y)^2)((x+y)^2-z^2)\right]^{\alpha-\nicefrac12}\;dz\\
      & =\frac{C_\Gamma}{xy}\int_{\abs{x-y}}^1z\Big[\underbrace{(z^2-(x-y)^2)}_{\leq1}\underbrace{\left(\frac{(x+y)^2-z^2}{x^2y^2}\right)}_{\leq\frac{(x+y)^2-(x-y)^2}{x^2y^2}=\frac4{xy}}\Big]^{\alpha-\nicefrac12}\;dz\\
      & \leq\frac{C_\Gamma}{xy}\int_{\abs{x-y}}^1z\left(\frac 4{xy}\right)^{\alpha-\nicefrac12}\;dz \leq \frac{C_\Gamma\cdot 4^{\alpha-\nicefrac12}}{(xy)^{\alpha+\nicefrac12}}.\\
\intertext{The last inequality is due to $1-\abs{x-y}\leq 1$ and the standard estimate. By using $k-1\leq x < k+2$ and $k\leq y < k+1$, we get}
      & \leq\frac{C_\Gamma4^{\alpha-\nicefrac12}}{(k-1)^{\alpha+\nicefrac12}k^{\alpha+\nicefrac12}}\leq\frac{C_\Gamma4^{\alpha-\nicefrac12}\left(\frac{k}{k-1}\right)^{\alpha+\nicefrac12}}{k^{\alpha+\nicefrac12}k^{\alpha+\nicefrac12}}\\
      & \leq\frac{C_\Gamma4^{\alpha-\nicefrac12}2^{\alpha+\nicefrac12}}{k^{\alpha+\nicefrac12}k^{\alpha+\nicefrac12}}\leq\frac{C_\Gamma4^{\alpha-\nicefrac12}2^{\alpha+\nicefrac12}}{k^{2\alpha+1}}\\
      & = C\cdot\frac1{k^{2\alpha+1}}.
\end{align}
    We note that
    \begin{equation}
      \omega_{k+2}=\int_{k+1}^{k+2}z^{2\alpha+1}\;dz\leq(k+2)^{2\alpha+1}\leq 2^{2\alpha+1}k^{2\alpha+1}.
    \end{equation}
    Hence we get
    \begin{equation}
      \tau_y\eins_{[0,1]}(x)\leq \frac {C'}{\omega_{k+2}}\leq\frac {C'}{\omega_{k+1}}\leq\frac {C'}{\omega_{k}}\quad\text{for }y\in I_k,\, k\geq2.
    \end{equation}
    Altogether we obtain
    \begin{align*}
      \int_{I_j}\abs f^p\tau_y\eins_{[0,1]}\;d\omega_\alpha&\leq C'\cdot\int_{I_j}\frac1{\omega_j}\abs f^p\;d\omega_\alpha\quad  (j=k,k+1,k+2)\\
      \Rightarrow\left(\int\abs f^p\tau_y\eins_{[0,1]}\;d\omega_\alpha\right)^{\nicefrac1p}&\leq C'^{\nicefrac1p}\cdot\sum_{j=k}^{k+2}\left(\int_{I_j}\frac1{\omega_j}\abs f^p\;d\omega_\alpha\right)^{\nicefrac1p}\\
      &\leq3\cdot C'^{\nicefrac1p}\cdot\norm f_{p,\infty}
    \end{align*}
Choosing $C''$ as the maximum of the constants in i) - iii) concludes the proof.
  \end{enumerate}
\end{bew}

\subsection{Uniform boundedness of the translation operator on $(L^p,\ell^q)(\R_+,\conv_\alpha)$}
\begin{pro}[The case $(L^\infty,\ell^1)(\R_+,\conv_\alpha)$]\label{pro:equiv}
For $f\in(L^\infty,\ell^1)(\R,\conv_\alpha)$ and $y\in\R_+$ it holds that
\begin{equation}
  \norm{\tau_y f}_{\infty,1}\leq C\cdot\norm f_{\infty,1} 
\end{equation}
with a constant $C$ independent of $y$.
\end{pro}
\begin{bew}
Like in \cite{leinpaper} it is enough to consider only functions $f_n\coloneqq\eins_{I_n}=\eins_{[n-1,n)}$. That means we show only that
\begin{equation}
  \norm{\tau_yf_n}_{\infty,1}\leq C\cdot\norm{f_n}_{\infty,1}=C\cdot\omega_n
\end{equation}
with $C$ independent of $n$ and $y$. For the sake of completeness, we repeat the proof for the correctness of our constraint.\\
Let  $c_n\coloneqq\norm{P_n f}_\infty$ ($P_n$ denotes the restriction to the interval $I_n$) and let $g=\sum_nc_nf_n$. 
\begin{align}
  \norm{\tau_yf}_{\infty,1}&\leq\norm{\tau_yg}_{\infty,1}\leq\sum_nc_n\norm{\tau_yf_n}_{\infty,1}\\
  &\leq\sum_nc_nC\norm{f_n}_{\infty,1}=C\norm f_{\infty,1}.
\end{align}
Fix $y$ and $n$. We denote a $k\in\Z_+$ as \emph{exceptional} if $k=1$ or if there exists $x\in I_k$ so that $\abs{x-y}$ or $x+y$ lies in $I_n$. The set of all exceptional indices will be denoted as $E$. An index $k\in\Z_+$ which is not exceptional will be called \emph{generic}. $G\coloneqq\Z_+\setminus E$.\\
For $k\in G$ the intersection of $[\abs{x-y},x+y]$ and $I_n$ is either empty or all of $I_n$ for all $x\in I_k$. Then $\tau_yf_n$ either vanishes on all of $I_k$ or is in the form of
\begin{equation}
  \tau_yf_n(x)=\frac{C_\Gamma}{(xy)^{2\alpha}}\int_{n-1}^nz\left[\left(z^2-(x-y)^2\right)\left((x+y)^2-z^2\right)\right]^{\alpha-\nicefrac12}\;dz,\quad\text{for }x\in I_k.
\end{equation}
So it is easy to see that we can claim the following statements about $x\in \{\tau_y f_n >0\}$ if $x\in I_k$ and $k$ generic.
\begin{gather}
  \abs{x-y}<n-1<n\leq x+y\\
\Rightarrow \quad x-y<n-1 \;\text{ and }\; y-x<n-1\;\text{ and }\; n\leq x+y\\
\Rightarrow \quad x<n+y-1 \;\text{ and }\; y-n+1<x\;\text{ and }\; n-y\leq x.
\end{gather}
Taken together it holds that $\supp\tau_y f_n\cap\left(\bigcup_{k\in G} I_k\right)\subset[\max\{0,y-n+1,n-y\},n+y-1]$. Equating the two terms ($\neq0$) of the lower bound yields $y=n-\nicefrac12$. Thus we can distinguish two cases: 

\begin{enumerate}
\item $y\leq n-\nicefrac12\Rightarrow n-y\leq x < n+y-1$ \label{fall1}
\item $y\geq n-\nicefrac12\Rightarrow y-n+1 < x < n+y-1$.\label{fall2}
\end{enumerate}
Using  $\omega_k\leq k^{2\alpha+1}$ and\\ $\frac1{x^{2\alpha}}\leq \frac1{(k-1)^{2\alpha}}\leq \frac{\left(\frac{k}{k-1}\right)^{2\alpha}}{k^{2\alpha}}\leq C\cdot\frac 1 {k^{2\alpha}} $ we get for $x\in I_k,\,k\geq2$
\begin{align*}
\MoveEqLeft\frac1{\omega_n}\sum_{k\in G}\omega_k\norm{P_k(\tau_y f_n)}_\infty\\
= {} & C_\Gamma\cdot\sum_{k\in G}\frac{\omega_k}{\omega_n}\sup_{x\in I_k}\int_{n-1}^n\frac z{(xy)^{2\alpha}}\left[\left(z^2-(x-y)^2\right)\left((x+y)^2-z^2\right)\right]^{\alpha-\nicefrac12}\;dz\\
\leq {} & C'\cdot\sum_{k\in G}\frac{k^{2\alpha+1}}{\omega_ny^{2\alpha} k^{2\alpha}}\sup_{x\in I_k}\int_{n-1}^nz\left[\left(z^2-(x-y)^2\right)\left((x+y)^2-z^2\right)\right]^{\alpha-\nicefrac12}\;dz\\
\leq {} & C'\cdot\sum_{k\in G}\frac{k}{\omega_ny^{2\alpha}}\sup_{x\in I_k}\int\limits_{n-1}^nz\big[\underbrace{(z-x+y)(z+x-y)(x+y-z)(x+y+z)}_{\eqqcolon\circledast}\big]^{\alpha-\nicefrac12}\;dz.\label{transl1}\tag{\ensuremath{\ast}}\\
\intertext{We consider $n$ and $y$ as in the first case and substitute, according to the sign, $z$ by $n$ resp. $n-1$ and $x$ by $n-y$ resp. $n+y-1$. It follows that}
\eqref{transl1} \leq {} & C'\cdot\sum_{k\in G}\frac k{\omega_n y^{2\alpha}}\sup_{x\in I_k} n\big[(n-n+y+y)(n+n+y-1-y)\\
  & \qquad\qquad (n+y-1-n+1+y)(n+y-1+y+n)\big]^{\alpha-\nicefrac12}\\
\leq {} & C'\cdot\sum_{k\in G}\frac k{\omega_n y^{2\alpha}}n\left[2y(2n-1)2y(2n+2n-2)\right]^{\alpha-\nicefrac12}\qquad{\scriptstyle(\text{because }2y\leq 2n-1)}\\
\leq {} & C''\cdot\sum_{k\in G}\frac{kny^{2\alpha-1}n^{2\alpha-1}}{n^{2\alpha+1}y^{2\alpha}}=C''\cdot\sum_{k\in G}\frac{k}{ny}\\
\leq {} & C''\cdot\sum_{k\in G}\frac2y\leq 4C''.
\end{align*}
In the last row it was used that $k \leq n+y-1 \leq 2n$ and that the sum consists of fewer than $2y$ terms: $n+y-1 - (n-y)=2y-1$.\\ 
Now let $y$ and $n$ be as in the second case. Then we can estimate $\circledast$ analogously.
\begin{align*}
  \circledast & \leq \big[(n-y+n-1+y)(n+n+y-1-y)\\
  & \qquad\qquad(n+y-1+y-n+1)(n+y-1+y+n)\big]^{\alpha-\nicefrac12}\\
  & \leq \left[(2n-1)(2n-1)2y(2y+2y)\right]^{\alpha-\nicefrac12}\qquad {\scriptstyle(\text{because }2n-1\leq 2y)}\\
  & \leq C\cdot n^{2\alpha-1}y^{2\alpha-1}\\
\Rightarrow\eqref{transl1}&\leq C''\cdot\sum_{k\in G}\frac{kn^{2\alpha}y^{2\alpha-1}}{y^{2\alpha}n^{2\alpha+1}}=C''\cdot\sum_{k\in G}\frac{k}{ny}\\
& \leq C''\cdot \sum_{k\in G}\frac2n\leq 4C''.
\end{align*}
Similar to the first case we used here that $k\leq n+y-1\leq 2y$ and that the sum consists of fewer than $n+y-1-(y-n+1)=2n-2$ terms.

If $k$ is now an exceptional index, we have to estimate $\omega_k\norm{P_k(\tau_y f_n)}_\infty$. For exceptional indices we have by definition either $x+y\in I_n$ for $x\in I_k$, i.e. $y+I_k$ intersects the interval $I_n$, or $I_n-y$ intersects $I_k$. There are at most two such indices $k\in E$. The other cases of exceptional indices derive from the cases where  $I_n+y$ or $y-I_n$ intersect the interval $I_k$, or where $k=1$. Each of these cases, with the exception of $k=1$, yields at most $2$ exceptional indices. Thus there are at most $7$. By looking closer one can see that there are in fact only $5$ of them.\\

If $k\leq 3n$ then one has, with the use of
\begin{align*}
  \tau_y f_n(x) & = \int_{\abs{x-y}}^{x+y} f_n(x) K(x,y,z)\;\omega_\alpha(dz)\leq\int_{\abs{x-y}}^{x+y} K(x,y,z)\;\omega_\alpha(dz)=1\\
\shortintertext{and}
\omega_{3n} & = \int_{3n-1}^{3n}z^{2\alpha+1}\;dz =\frac 1{2\alpha+2}\left((3n)^{2\alpha+2}-(3n-1)^{2\alpha+2}\right) \\
    & = \frac{3^{2\alpha+2}}{2\alpha+2}\left(n^{2\alpha+2}-(n-\frac13)^{2\alpha+2}\right)\\
    & \leq \frac {3^{2\alpha+2}}{2\alpha+2}\left(n^{2\alpha+2}-(n-1)^{2\alpha+2}\right)=3^{2\alpha+2}\int_{n-1}^nz^{2\alpha+1}\;dz \\
    & = 3^{2\alpha+2}\omega_n\;,
\end{align*}
the following estimate:
\begin{equation}
\omega_k\norm{P_k(\tau_yf_n)}_\infty  \leq \omega_k\leq\omega_{3n}\leq 3^{2\alpha+2}\omega_n=3^{2\alpha+2}\norm{f_n}_{\infty,1}.
\end{equation}
If $k$ is exceptional and $k>3n$, then one of the intervals $y\pm I_n$ must intersect $I_k$. The smallest value for $y$ must therefore satisfy $y+n=k-1$. That implies that
\begin{equation}
  y+\frac13k>k-1\Rightarrow y>\frac23k-1>\frac13k,
\end{equation}
because $k>3$. Particularly we have $y>\frac13x$ for all $x\in I_k$ in these cases. Now we can find an upper bound for the remaining exceptional indices.
\begin{align*}
  \tau_yf_n(x)&=\frac{C_\Gamma}{(xy)^{2\alpha}}\int_{\abs{x-y}}^nz\Big[\underbrace{\left(z^2-(x-y)^2\right)}_{\leq n^2-(x-y)^2}\underbrace{\left((x+y)^2-z^2\right)}_{\mathclap{\leq (x+y)^2-\abs{x-y}^2=4xy}}\Big]^{\alpha-\nicefrac12}\;dz\\
  & \leq\frac{C_\Gamma}{(xy)^{2\alpha}}\underbrace{(n-\abs{x-y})}_{\leq 1} \cdot n\cdot \left[\left(n^2-(x-y)^2\right)4xy\right]^{\alpha-\nicefrac12}\\
  & \leq\frac{C_\Gamma4^{\alpha-\nicefrac12}}{(xy)^{\alpha+\nicefrac12}}n\Big[\underbrace{n^2-(x-y)^2}_{\mathclap{\leq n^2-(n-1)^2=2n-1}}\Big]^{\alpha-\nicefrac12}\\
  & \leq\frac{C_\Gamma4^{\alpha-\nicefrac12}}{(\frac13x*x)^{\alpha+\nicefrac12}}n\underbrace{(2n-1)}_{\leq 2n}{}^{\alpha-\nicefrac12}\\
  & \leq\frac{C_\Gamma4^{\alpha-\nicefrac12}3^{\alpha+\nicefrac12}}{x^{2\alpha+1}}2^{\alpha-\nicefrac12}n^{\alpha+\nicefrac12}\\
  & \leq C\frac{n^{\alpha+\frac12}}{(k-1)^{2\alpha+1}}.
\end{align*}
Hence we can estimate the remaining terms of the norm.
\begin{equation*}
  \omega_k\norm{P_k(\tau_yf_n)}_\infty\leq C\cdot\frac{k^{2\alpha+1}n^{\alpha+\nicefrac12}}{(k-1)^{2\alpha+1}}\leq C\cdot \left(\frac43\right)^{\mathclap{\quad2\alpha+1}}\cdot n^{\alpha+\nicefrac12}\leq C'\cdot\omega_n=C'\cdot\norm{f_n}_{\infty,1}
\end{equation*}
\end{bew}
Using duality and complex interpolation as in \cite{leinpaper}, the boundedness of translation on $(L^\infty,\ell^1)(\R_+,\conv_\alpha)$ can be extended to the general amalgam spaces $(L^p,\ell^q)(\R_+,\conv_\alpha)$. 
\subsection{Fourier transformation on $(L^p,\ell^q)(\R_+,\conv_\alpha)$}

We have yet to prove the following theorem:
\begin{satz}[Hausdorff-Young theorem for amalgams]
  For $f\in(L^p,\ell^q)(\R_+,\conv_\alpha)$ with $1\leq p,q\leq 2$ holds, that $\hat f\in(L^{q'},\ell^{p'})(\R_+,\conv_\alpha)$.
\end{satz}
The special case $p=q=2$ is already known (see Theorem 2.2.22 in \cite{bloomheyer}). We will consider the other extreme cases $p,q\in\{1,2\}$ first.
Therefore we will have to show that $\norm{\widehat{\eins_{I_1}}}_{\infty,2}<\infty$.
\begin{align}
  \widehat{\eins_{I_1}}(\lambda)&= \int_{\R_+}\eins_{I_1}(x)\chi_\lambda(x)\,d\omega_\alpha(x)\\
  &=\int_0^1\sum_{k=0}^\infty\frac{(-1)^k\Gamma(\alpha+1)}{2^{2k}k!\Gamma(\alpha+k+1)}(\lambda x)^{2k}\,d\omega_\alpha(x)\\
  &=\sum_{k=0}^\infty\frac{(-1)^k\Gamma(\alpha+1)\lambda^{2k}}{2^{2k}k!\Gamma(\alpha+k+1)}\int_0^1x^{2k}x^{2\alpha+1}\,dx\\
  &=\sum_{k=0}^\infty\frac{(-1)^k\Gamma(\alpha+1)\lambda^{2k}}{2^{2k}k!\Gamma(\alpha+k+1)}\frac1{2k+2\alpha+2}\\
  &=\sum_{k=0}^\infty\frac{(-1)^k\Gamma(\alpha+1)\lambda^{2k}}{2^{2k+1}k!\Gamma(\alpha+k+2)}\\
  &=\frac{\Gamma(\alpha+1) 2^\alpha}{\lambda^{\alpha+1}}\sum_{k=0}^\infty\frac{(-1)^k\lambda^{2k+\alpha+1}}{2^{2k+\alpha+1}k!\Gamma(\alpha+k+2)}\\
  &=\frac{\Gamma(\alpha+1) 2^\alpha}{\lambda^{\alpha+1}}J_{\alpha+1}(\lambda).
\end{align}
Here $J_{\alpha+1}$ denotes the Bessel function of degree $\alpha+1$. According to 9.2.1 in \cite{Handbook} it holds for $\lambda\to\infty$ that
\begin{equation}
  J_{\alpha+1}(\lambda)\approx\sqrt{\frac2{\pi\lambda}}\cos(\lambda-\tfrac12(\alpha+1)\pi-\tfrac14\pi).
\end{equation}
Altogether for large enough $\lambda$ we have
\begin{equation}
\abs*{\widehat{\eins_{I_1}}(\lambda)} \approx \abs*{\frac{\Gamma(\alpha+1) 2^\alpha}{\lambda^{\alpha+1}}\sqrt{\frac2{\pi\lambda}}\cos(\lambda-\tfrac12(\alpha+1)\pi-\tfrac14\pi)} \lessapprox C\cdot \lambda^{-\alpha-\nicefrac32}.
\end{equation}

Now we can show that $\norm{\widehat{\eins_{I_1}}}_{\infty,q}<\infty\gdw q>2\frac{\alpha+1}{\alpha+\nicefrac32}$. Particularly for $q\geq2$
\begin{align}
  \norm{\widehat{\eins_{I_1}}}_{\infty,q} & =\left(\sum_{k=0}^\infty\omega_k\sup_{\lambda\in I_k}\abs{\widehat{\eins_{I_1}}(\lambda)}^q\right)^{\nicefrac1q}\\
  & \leq C\cdot\left(C'+\sum_{k=N_\alpha}^\infty k^{2\alpha+1}k^{q(-\alpha-\nicefrac32)}\right)^{\nicefrac1q}\\
  & = C''+C\cdot\left(\sum_{k=N_\alpha}^\infty k^{2\alpha+1-q(\alpha+\nicefrac32)}\right)^{\nicefrac1q},
\end{align}
where $N_\alpha$ and $C$ be chosen such that $\abs*{\widehat{\eins_{I_1}}(\lambda)} \leq C\cdot\lambda^{-\alpha-\nicefrac32}$ for $\lambda\geq N_\alpha$. The series converges iff $2\alpha+1-q(\alpha+\nicefrac32)<-1\gdw q>2\frac{\alpha+1}{\alpha+\nicefrac32}$. One can easily see that $\widehat{\eins_{I_1}}$ does not lie in $(L^p,\ell^q)(\R_+,\conv_\alpha)$ if $q\leq 2\frac{\alpha+1}{\alpha+\nicefrac32}$ and for all $p$.
Now we proceed again as in \cite{leinpaper}. Let $g_1=\frac1{\omega_1}\cdot\eins_{I_1}\conv_\alpha\eins_{[0,2]}$ and $g_n=\frac1{\omega_1}\cdot\eins_{I_1}\conv_\alpha\eins_{[n-2,n+1]}$ for $n\geq2$. We show that $g_n(x)=1\;\forall x\in I_n$. For this let $x\in I_1$. 
\begin{align}
  g_1(x)&=\frac1{\omega_1}\int_{\R_+}\eins_{I_1}(y)\tau_x\eins_{[0,2]}(y)\,d\omega_\alpha(y)\\
  &=\frac1{\omega_1}\int_0^1\int_{\abs{x-y}}^{x+y}\eins_{[0,2]}(z)K(x,y,z)\,d\omega_\alpha(z)d\omega_\alpha(y)\\
  &=\frac1{\omega_1}\int_0^1\int_{[\abs{x-y},x+y]\cap[0,2]}K(x,y,z)\,d\omega_\alpha(z)d\omega_\alpha(y)\\
  &=\frac1{\omega_1}\int_0^11\,d\omega_\alpha(y)=1,
\end{align}
because for $x,y\in[0,1]\Rightarrow [\abs{x-y},x+y]\subset[0,2]\Rightarrow[\abs{x-y},x+y]\cap[0,2]=[\abs{x-y},x+y]$ and thus the inner integral yields the value $1$. Entirely analog it holds for $x\in I_n$, that
\begin{align}
  g_n(x)&=\frac1{\omega_1}\int_{\R_+}\eins_{I_1}(y)\tau_x\eins_{[n-2,n+1]}(y)\,d\omega_\alpha(y)\\
  &=\frac1{\omega_1}\int_0^1\int_{[\abs{x-y},x+y]\cap[n-2,n+1]}K(x,y,z)\,d\omega_\alpha(z)d\omega_\alpha(y)\\
  &=\frac1{\omega_1}\int_0^11\,d\omega_\alpha(y)=1,
\end{align}
because for $x\in[n-1,n]$ and $y\in[0,1]\Rightarrow [\abs{x-y},x+y]\subset[n-2,n+1]\Rightarrow[\abs{x-y},x+y]\cap[n-2,n+1]=[\abs{x-y},x+y]$.

The rest of the proof runs exactly as in \cite{leinpaper}, with $3$ replaced by $\frac{1}{\omega_1}$ and $\conv_{\nicefrac12}$ replaced by $\conv_\alpha$. Note that Young's inequality can be obtained as in \cite{leinpaper} as well. 
\hfill$\square$

\newpage

\section{Equivalence on compact/discrete hypergroups}

As an addendum, if $K$ is a compact or discrete Hypergroup, let us note the equivalence of the continuous amalgam norm (here denoted as $\norm{\cdot}^*_{p,\infty}$) with the discrete amalgam norm. \\
We first describe the discrete case. The discrete amalgam norm on a discrete hypergroup is defined as one would expect as
\begin{align}
  \norm f_{p,q} & = \left(\sum_{k\in K}\omega(\{k\})\left(\frac1{\omega(\{k\})}\abs{f(k)}^p\omega(\{k\})\right)^{\nicefrac q p}\right)^{\nicefrac1q}\\
  & =\left(\sum_{k\in K}\omega(\{k\})\abs{f(k)}^q\right)^{\nicefrac1q}=\norm f_q.
\end{align}
For the equivalence of the norms we compute
\begin{align}
  \norm f_{p,\infty}^* & =\sup_{n\in K}\left(\sum_{k\in K}\tau_n\eins_{\{0\}}(k)\abs{f(k)}^p\omega(\{k\})\right)^{\nicefrac 1 p}\\
  & = \sup_{n\in K}\left(\sum_{k\in K}\eins_{\{0\}}(n\conv k)\abs{f(k)}^p\omega(\{k\})\right)^{\nicefrac 1 p}\\
\intertext{so, using $0\in\supp\eps_n\conv\eps_k\gdw k=n\inv$, we get }
  & = \sup_{n\in K}\Big(\eps_{n\inv}\conv\eps_{n}(\{0\})\abs{f(n)}^p\omega(\{n\})\Big)^{\nicefrac 1 p}\\
\intertext{and finally, with Theorem 1.3.26 in \cite{bloomheyer}:  $\omega(\{n\})=\eps_{n\inv}\conv\eps_{n}(\{0\})^{-1}$, we have}
  & = \sup_{n\in K}\abs{f(n)}=\norm f_\infty.
\end{align}
On the other hand, in the compact case, the discrete amalgam norm shrinks to
\begin{equation}
 \left(\omega(K)\left(\frac1{\omega(K)}\int_K\abs f^p\;d\omega\right)^{\nicefrac q p}\right)^{\nicefrac1 q}=\left(\int_K\abs f^p\;d\omega\right)^{\nicefrac 1 p}=\norm f_p.
\end{equation}
The equivalence now follows:
\begin{align}
  \norm f_{p,\infty}^* & =  \sup_{y\in K}\left(\int\tau_y\eins_{K}\abs f^p\;d\omega\right)^{\nicefrac1p}\\
    & = \sup_{y\in K}\left(\int\int_{K}\eins_{K}\;d(\eps_y\conv\eps_x) \abs f^p(x)\;\omega(dx)\right)^{\nicefrac1p}\\
    & = \sup_{y\in K}\left(\int\abs f^p\;d\omega\right)^{\nicefrac1p}=\left(\int\abs f^p\;d\omega\right)^{\nicefrac1p}=\norm f_p.
\end{align}
So, in both cases, we obtained not only equivalence, but equality of norms.


\end{document}